\documentclass[11pt]{amsart}
\usepackage{amsmath,amsthm,epsfig,amssymb,graphicx}
\title{Introduction to the Prisoners Versus Guards Game}
\author{Timothy Howard, Eugen J. Ionascu, and David Woolbright }
\curraddr{Math Department\\Columbus State University\\4225 University Avenue\\
Columbus, GA 31907}

\email{thoward@colstate.edu;}
\email{ionascu\_eugen@colstate.edu;}
\email{woolbright\_david@colstate.edu;}

\date{January 21, 2007}
\begin{document}
\def\sms{\small\scshape}
\baselineskip18pt
\def\RR{{\rm I}\!{\rm R}}
\def\fp#1{(#1)}
\newtheorem{theorem}{\hspace{\parindent}
T{\scriptsize HEOREM}}[section]
\newtheorem{proposition}[theorem]
{\hspace{\parindent }P{\scriptsize ROPOSITION}}
\newtheorem{corollary}[theorem]
{\hspace{\parindent }C{\scriptsize OROLLARY}}
\newtheorem{lemma}[theorem]
{\hspace{\parindent }L{\scriptsize EMMA}}
\newtheorem{definition}[theorem]
{\hspace{\parindent }D{\scriptsize EFINITION}}
\newtheorem{problem}[theorem]
{\hspace{\parindent }P{\scriptsize ROBLEM}}
\newtheorem{conjecture}[theorem]
{\hspace{\parindent }C{\scriptsize ONJECTURE}}
\newtheorem{example}[theorem]
{\hspace{\parindent }E{\scriptsize XAMPLE}}
\newtheorem{remark}[theorem]
{\hspace{\parindent }R{\scriptsize EMARK}}
\renewcommand{\thetheorem}{\arabic{section}.\arabic{theorem}}
\newcommand{\du}{\stackrel{.}{\bigcup}}
\renewcommand{\theenumi}{(\roman{enumi})}
\renewcommand{\labelenumi}{\theenumi}
\def\RR{{\rm I}\!{\rm R}}
\def\LL{{\rm I}\!{\rm L}}
\def\NN{{\rm I}\!{\rm N}}
\def\MM{{\rm I}\!{\rm M}}
\def\QQ{{\rm I}\!\!\!{\rm Q}}
\def\ZZ{{\rm Z}\!\!{\rm Z}}
\def\CC{{\rm I}\!\!\!{\rm C}}
\newcommand{\Q}{{\mathbb Q}}
\newcommand{\Z}{{\mathbb Z}}
\newcommand{\N}{{\mathbb N}}
\newcommand{\C}{{\mathbb C}}
\newcommand{\R}{{\mathbb R}}
\newcommand{\F}{{\mathbb F}}
\newcommand{\K}{{\mathbb K}}
\newcommand{\D}{{\mathbb D}}
\def\vp{\varepsilon}
\def\phi{\varphi}
\def\ra{\rightarrow}
\def\sd{\bigtriangledown}
\def\ac{\mathaccent94}
\def\wi{\sim}
\def\wt{\widetilde}
\def\bb#1{{\Bbb#1}}
\def\bs{\backslash}
\def\cal{\mathcal}
\def\ca#1{{\cal#1}}
\def\Bbb#1{\bf#1}
\def\blacksquare{{\ \vrule height7pt width7pt depth0pt}}
\def\bsq{\blacksquare}
\def\proof{\hspace{\parindent}{P{\scriptsize ROOF}}}
\def\pofthe{P{\scriptsize ROOF OF}
T{\scriptsize HEOREM}\  }
\def\pofle{\hspace{\parindent}P{\scriptsize ROOF OF}
L{\scriptsize EMMA}\  }
\def\pofcor{\hspace{\parindent}P{\scriptsize ROOF OF}
C{\scriptsize ROLLARY}\  }
\def\pofpro{\hspace{\parindent}P{\scriptsize ROOF OF}
P{\scriptsize ROPOSITION}\  }
\def\n{\noindent}
\def\eproof{$\hfill\bsq$\par}
\def\ds{\displaystyle}
\def\du{\overset{\text {\bf .}}{\cup}}
\def\Du{\overset{\text {\bf .}}{\bigcup}}
\def\b{$\blacklozenge$}

\maketitle

\section{\label{intro}Introduction}

Suppose that you are competing in a two-player game in which you and
your opponent attempt to pack as many ``prisoners'' as possible on
the squares of an $n\times n$ checkerboard; each prisoner has to be
``protected'' by an appropriate number of guards.  Initially, the
board is covered entirely with guards.  The players -- designated as
``red'' and ``blue'' -- take turns adjusting the board configuration
using one of the following rules in each turn:
\begin{itemize}
\item [I.]  Replace one guard with a prisoner of the player's
color.
\item [II.]  Replace one prisoner of either color with a guard
and replace two other guards with prisoners of the player's color.
\end{itemize}

\begin{figure}[h]
\centering
\includegraphics[scale=0.25]{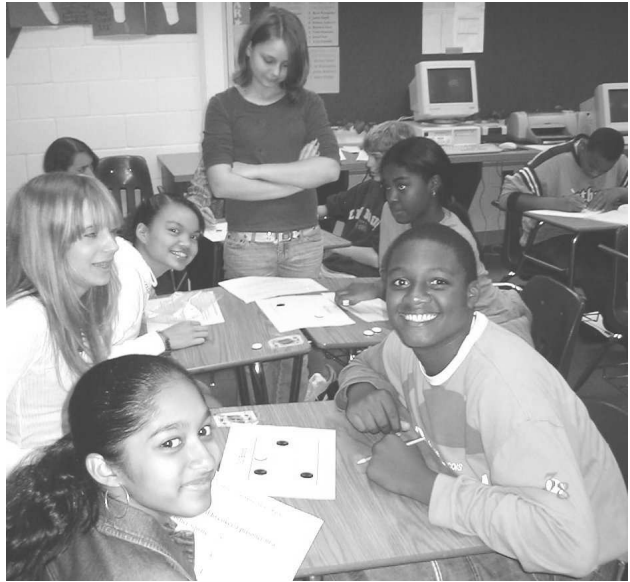}
\caption{Blackmon Road Middle School Students Play PvG}
\label{blackmon-rd-players}
\end{figure}

\noindent That is, each player takes a turn increasing the total
number of prisoners by one.  We require that, at every stage of the
game, {\it each prisoner lies adjacent to at least as many guards as
the number of the other prisoners adjacent it}.  The squares {\it
adjacent} to a given square are those squares, situated directly
above, below, to the left, to the right, or diagonal to the square
in question. An arrangement of prisoners and guards that satisfies
this requirement and has exactly one occupant per square is called a
{\it valid board}. The game ends when neither player can further
adjust the board using rules I and II while maintaining a valid
board.  The player whose color represents more prisoners is the
winner.  This is the game of Prisoners and Guards -- a game that can
be played and analyzed without extensive knowledge of mathematics.
Figure~\ref{blackmon-rd-players} depicts students from Blackmon Road
Middle School, in Columbus, Georgia, trying their hands at the game.
 We invite the reader to play the game online by running the Java
Applet available at http://csc.colstate.edu/woolbright/.

The guards in this game are related to the half domination set in
the king's graph as introduced in in a paper by Hoffman, Laskar, and
Markus (see \cite{dhlm}). Similar domination problems are studied
have been studied by Bode, Harborth, and Harporth (see \cite{bhh})
and by Watkins, Ricci, and McVeigh (see \cite{wrm}). The Prisoners
and Guards game originated as a puzzle created by the third author
with a focus on minimizing the size of the dominating set (the
guards).

In the two-player game, one fundamental question that naturally
arises is ``How do we decide when the game is over?''  The short
answer is that the game is over when the board configuration has
reached a maximal state.  A valid board to which no adjustments can
be made to increase the total number of prisoners is called a {\it
maximal} board. One can also define a {\it maxi\underline{mum}}
configuration as being a maximal arrangement of prisoners and guards
that has the greatest number of prisoners of all valid boards.  For
$n\in \left\{1,2,3\right\}$ all maximal boards are also maximum
configurations.  We will see example $4\times 4$ boards that are
maximal but do not contain maximum configurations.  Let $P(n)$
denote the number of prisoners in a maximum configuration.
Characterizing the sequence $\left\{P(n)\right\}_{n=1}^\infty$ will
help us determine when to end the game.  It also proves an
interesting avenue for exploration on its own.

Since the lone square on a $1\times 1$ board has no adjacent
squares, we can place a prisoner in it and be sure that there are at
least as many guards as prisoners lying in all adjacent squares --
none. Therefore, we have $P(1)=1$.  By exhaustively checking all
sixteen $2\times 2$ cases, we find eleven valid boards, each having
zero, one, or two prisoners.  Thus, $P(2)=2$.  We analyze the cases
$n=3$, $4$, $5$, and $6$ in sections \ref{3x3and4x4cases} and
\ref{5x5and6x6cases}.  Exact values for $P(n)$, $n\in
\left\{7,8,9,10,11\right\}$, can be found in a paper by Ionascu,
Pritikin, and Wright (see \cite{ipw}), who employ linear programming
techniques in the study of $P(n)$. In Section~\ref{deficiency} we
also obtain an upper bound on $P(n)$; this is about the best that we
can say for $n\ge 12$.

\section{\label{3x3and4x4cases}The game analysis for $n=3$ and $n=4$}

\begin{figure}[ht]
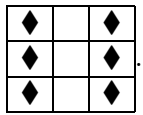

\begin{tabular}{|r|r|r|}
\hline
\b & \ \  & \b \\
\hline
\b & \ \  & \b \\
\hline
\b &\ \  & \b \\
\hline
\end{tabular}.
\caption{Maximum 3x3 Board}
\label{3x3maximumboard}
\end{figure}

Playing Prisoners and Guards on a $1\times 1$ board or on a $2\times
2$ board is not all that interesting.  When we increase the board
size slightly and consider the game on a $3\times 3$ board, strategy
becomes more of a factor.  We will see that the arrangement in
Figure~\ref{3x3maximumboard} is a maximum configuration, as we
establish in Theorem~\ref{3x3maxtheorem}.  We use the diamond to
represent prisoners and blank squares represent guards.

\noindent In fact, this is the only maximal arrangement (up to a
rotation).  Let us observe that a maximal board permits no
adjustments using either Rule I or Rule II.  First we consider
arrangements that are maximal with respect to Rule I (i.e. one
cannot simply add more prisoners in the existing configuration).

\begin{figure}[ht]
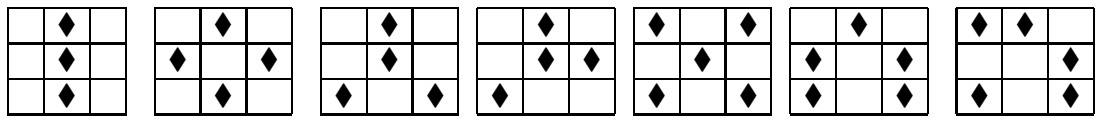


\begin{tabular}{|r|r|r|}
\hline
\ \ & \b \  & \ \ \\
\hline
\ \ & \b \  & \ \ \\
\hline
\ \ &\b \  & \ \ \\
\hline
\end{tabular} \ \
\begin{tabular}{|r|r|r|}
\hline
\  & \b \  & \ \\
\hline
\b & \ \  & \b \\
\hline
\ &\b \  & \ \\
\hline
\end{tabular} \ \
\begin{tabular}{|r|r|r|}
\hline
\ \ & \b \  & \ \ \\
\hline
\ \ & \b \  & \ \ \\
\hline
\b \ &\ \  & \b \ \\
\hline
\end{tabular}\ \
\begin{tabular}{|r|r|r|}
\hline
\ \ & \b \  & \ \ \\
\hline
\ \ & \b \  & \b \ \\
\hline
\b \ &\ \  & \ \ \\
\hline
\end{tabular}\ \
\begin{tabular}{|r|r|r|}
\hline
\b & \ \  & \b \\
\hline
\ & \b \  & \ \\
\hline
\b &\ \  & \b \\
\hline
\end{tabular}\ \
\begin{tabular}{|r|r|r|}
\hline
\ & \b \  & \ \\
\hline
\b & \ \  & \b \\
\hline
\b &\ \  & \b \\
\hline
\end{tabular}\ \ \
\begin{tabular}{|r|r|r|}
\hline
\b \ & \b \  & \ \ \\
\hline
\ \ & \ \  & \b \ \\
\hline
\b \ &\ \  & \b \ \\
\hline
\end{tabular}
\caption{$3\times 3$ Boards that are Maximal w.r.t. Rule I}
\label{3x3maximalboardswrtruleI}
\end{figure}

Perhaps it is not difficult to convince oneself that any valid board
having zero, one, or two prisoners can be adjusted using Rule I;
after factoring out rotations and  reflections, there are ten unique
cases to check.  Therefore, a maximal board must have at least three
prisoners.  Figure~\ref{3x3maximalboardswrtruleI} depicts all valid
boards (up to rotations and reflections) that contain three, four,
or five prisoners and are maximal with respect to Rule I.

Each one of these configurations can be adjusted to match the
arrangement in Figure~\ref{3x3maximumboard} by using Rule II (and
one or more adjustments using Rules I and II in some of the cases).
It follows that a maximal $3\times 3$ board must contain at least
six prisoners.  We record this fact in the following lemma.

\begin{lemma} Any maximal $3\times 3$  board has at least six prisoners.
\label{3x3sixprisoners}\end{lemma}

With Figure~\ref{3x3maximumboard}, we see that a valid $3\times 3$
configuration can have six prisoners.  Does there exist a valid
configuration with more than six prisoners? Suppose that a $3\times
3$ board arrangement contains seven prisoners (and two guards).
Since there are four non-corner edge squares, a prisoner must occupy
at least one of them.  Since this prisoner lies adjacent to at most
two guards, the board cannot be valid.  These observations,
Lemma~\ref{3x3sixprisoners}, and the fact that
Figure~\ref{3x3maximumboard} depicts a valid configuration with six
prisoners lead us to the following conclusion.

\begin{theorem} A maximum $3\times 3$ board contains six prisoners, {\it i.e.}
$P(3)=6$. \label{3x3maxtheorem}\end{theorem}

As a matter of fact, the configuration shown in
Figure~\ref{3x3maximumboard} is the \underline{only} maximal
$3\times 3$ board (up to a rotation of the board).  From this we
learn that the second player has a good chance to win by using Rule
II all of the time. The first player may force a tie if she can lead
the board configuration in such a manner that will require her
opponent to use Rule I. In fact, this is manageable if she plays
into the pattern in Figure~\ref{3x3maximumboard}, forcing the blue
to use Rule I in the last step and so the number of prisoners of
each color will end up equal.

\begin{figure}[ht]
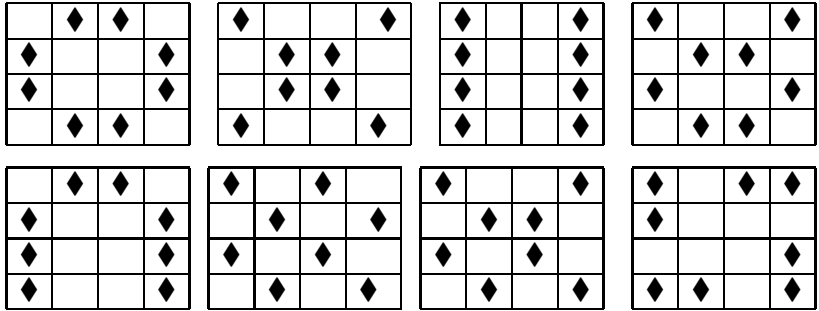


\begin{tabular}{|r|r|r|r|}
\hline
\    \       & \b   & \b & \   \\
\hline
\b   & \  \  & \  \  & \b \ \\
\hline
\b & \  \  & \   &    \b  \     \\
\hline
\ & \b \  & \b   & \  \  \   \\
\hline\end{tabular}\ \ \
\begin{tabular}{|r|r|r|r|}
\hline
\b    \       & \   & \ & \b  \\
\hline
\   & \b  \  & \b  \  & \ \ \\
\hline
\ & \b  \  & \b   &    \ \     \\
\hline
\b & \ \  & \   & \b  \  \   \\
\hline\end{tabular}\ \ \
\begin{tabular}{|r|r|r|r|}
\hline
\b   \ & \  \ & \ \  & \b   \\
\hline
\b  \ & \  \  & \  \  & \b \ \\
\hline
\b \ & \  \  & \  \  &    \b  \     \\
\hline
\b \ & \ \  & \ \  & \b  \    \\
\hline\end{tabular}\ \ \
\begin{tabular}{|r|r|r|r|}
\hline
\b   \ & \  \ & \ \  & \b   \\
\hline
\  \ & \b  \  & \b  \  & \ \ \\
\hline
\b \ & \  \  & \  \  &    \b  \     \\
\hline
\ \ & \b \  & \b \  & \  \    \\
\hline\end{tabular}

\vspace{0.1 in}
\begin{tabular}{|r|r|r|r|}
\hline
\   \ & \b  \ & \b \  & \   \\
\hline
\b  \ & \  \  & \  \  & \b \ \\
\hline
\b \ & \  \  & \  \  &    \b  \     \\
\hline
\b \ & \ \  & \ \  & \b  \    \\
\hline\end{tabular}\ \
\begin{tabular}{|r|r|r|r|}
\hline
\b   \       & \   & \b & \   \\
\hline
\  & \b  \  & \  \  & \b \ \\
\hline
\b & \  \  & \b   &    \  \     \\
\hline
\ & \b \  & \   & \b  \  \   \\
\hline\end{tabular}\ \
\begin{tabular}{|r|r|r|r|}
\hline
\b   \ & \  \ & \ \  & \b   \\
\hline
\  \ & \b  \  & \b  \  & \ \ \\
\hline
\b \ & \  \  & \b  \  &    \  \     \\
\hline
\ \ & \b \  & \ \  & \b  \    \\
\hline\end{tabular} \ \
\begin{tabular}{|r|r|r|r|}
\hline
\b   \ & \  \ & \b \  & \b   \\
\hline
\b  \ & \  \  & \  \  & \ \ \\
\hline
\ \ & \  \  & \  \  &    \b  \     \\
\hline
\b \ & \b \  & \ \  & \b  \    \\
\hline\end{tabular}
\caption{Some Maximal $4\times 4$ Boards}
\label{maximal4x4boards}
\end{figure}

We now turn our attention to $4\times 4$ boards.  This board size
proves interesting because there exist many maximal arrangements
that are not maximum configurations.  We include some known maximal
arrangements with eight prisoners in Figure~\ref{maximal4x4boards},
but there may be others.

\begin{figure}[ht]
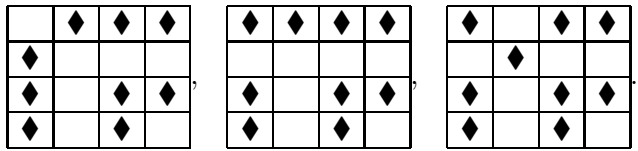

\begin{tabular}{|r|r|r|r|}
\hline
\ \       & \b   & \b & \b \\
\hline
\b   & \ \  & \ \  & \ \ \\
\hline
\b & \ \  & \b   & \b \\
\hline
\b & \ \  & \b   & \ \  \\
\hline
\end{tabular},\  \  \
\begin{tabular}{|r|r|r|r|}
\hline
\b  & \b   & \b & \b \\
\hline
\ \ & \ \  & \ \  & \ \  \\
\hline
\b & \ \   & \b  & \b \\
\hline
\b & \ \   & \b   & \ \  \\
\hline
\end{tabular},\  \   \
\begin{tabular}{|r|r|r|r|}
\hline
\b       & \ \   & \b & \b \\
\hline
\ \       & \b        & \ \  & \ \             \\
\hline
\b & \ \         & \b  & \b \\
\hline
\b & \ \         & \b  & \ \  \\
\hline
\end{tabular}.
\caption{Maximum $4\times 4$ Boards}
\label{4x4maximumboards}
\end{figure}

If we factor out rotations and reflections of the board, there are
three maximum arrangements as depicted in
Figure~\ref{4x4maximumboards}. We obtained these via an exhaustive
search of all $4\times 4$ valid boards and verified that there are
no other equivalence classes.

To prove something about maximum $4\times 4$ board configurations,
it helps to dissect the board and consider what can happen in the
vicinity of the corner squares.  Suppose that we have a $2\times 2$
block $C$ of squares situated in one corner of the board.  If the
corner square within $C$ does not contain a guard, then it contains
a prisoner.  If the latter is the case, then $C$ must contain at
least two guards.  Thus we have established the following fact.

\begin{lemma}\label{twobytwo}
If $C$ is a $2\times 2$ corner block within a valid board ($n \ge
2$), then it must contain at least one guard.
\end{lemma}

With this in mind, we are equipped to consider maximum $4\times 4$
boards by partitioning them into four $2\times 2$ corner blocks and
following through with the consequences.  This will lead us to the
conclusion of the next proposition.

\begin{proposition}\label{maximum4x4boards} $P(4)=9$.  That is, every maximum $4\times 4$
valid board has nine prisoners.
\end{proposition}

\begin{figure}[ht]
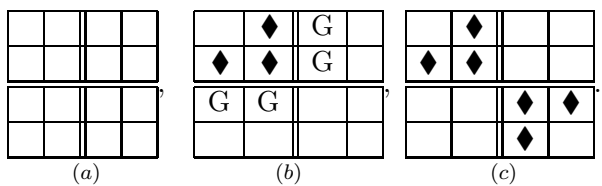

\[\underset{(a)} {\begin{tabular}{|r|r||r|r|}
\hline
\ \       & \    & \  & \  \\
\hline
\   & \ \  & \ \  & \ \ \\
\hline \hline
\ & \ \  & \   & \ \\
\hline
\ & \ \  & \   & \ \  \\
\hline
\end{tabular},} \ \ \
\underset{(b)}{
\begin{tabular}{|r|r||r|r|}
\hline
\  \       & \b    & G  & \  \\
\hline
\b   & \b \  & G   & \ \ \\
\hline \hline
G & G   & \   & \ \\
\hline
\ & \ \  & \   & \ \  \\
\hline
\end{tabular}}, \ \underset{(c)}{\begin{tabular}{|r|r||r|r|}
\hline
\  \       & \b    & \   & \  \\
\hline
\b   & \b \  & \   & \ \ \\
\hline \hline
\ & \ \  & \b   & \b \\
\hline
\ & \ \  & \b   & \ \  \\
\hline
\end{tabular}.}
\]
\caption{Block Partitions of a $4\times 4$ Board}
\label{blocksof4x4}
\end{figure}

\n \proof.   Since the configurations in
Figure~\ref{4x4maximumboards} are valid and each contains nine
prisoners, it follows that $P(4)\ge 9$. It is enough to show that
$P(4)\le 9$. We assume that there exist a valid board $B$ with ten
or more prisoners. By dropping prisoners if necessary, we can say
there are exactly ten.  We shall see that this leads to a
contradiction. We partition $B$ into four $2\times 2$ blocks as
indicated in Figure~\ref{blocksof4x4}(a).

Since by assumption the board contains ten prisoners, it follows
from Lemma~\ref{twobytwo} that at least two of these blocks must
contain three prisoners each. Without loss of generality, assume
that the upper left block is one of them. We see in the proof of
Lemma~\ref{twobytwo} that the lone guard must lie in square
$b_{11}$, as depicted in Figure~\ref{blocksof4x4}(b).

For the board to be valid, the non-corner edge prisoners in $b_{12}$
and $b_{21}$ must each lie adjacent to three guards. This is only
possible if guards lie in the squares $b_{13}$, $b_{23}$, $b_{31}$,
and $b_{32}$, as Figure~\ref{blocksof4x4}(b) indicates. As
previously noted, at least two of the   blocks must contain three
prisoners each. The only way to achieve this will be for the lower
right block to have three prisoners, with a guard in $b_{44}$ as
shown in Figure~\ref{blocksof4x4}(c).

Now we see that the prisoners situated in squares $b_{34}$ and
$b_{43}$ necessitate the presence of guards in squares $b_{24}$ and
$b_{42 }$. By placing prisoners in all squares not yet committed, we
will have a total of only eight prisoners on the board,
contradicting our assumption that the board has ten prisoners. Thus,
our assumption was invalid. \eproof

\section{\label{5x5and6x6cases}Analysis of the $5\times 5$ and $6\times 6$ Cases}

In our analysis of $5\times 5$ and $6\times 6$ board configurations,
we will partition the boards into $3\times 3$ blocks.  The following
lemma will help in the examination of these blocks.

\begin{lemma}\label{threebythree}
If $C$ is a $3\times 3$ corner block within a valid board ($n >
3$), then it must contain at least three guards. Moreover, if $C$
contains exactly three guards, then it must contain a prisoner
diagonally opposite (within $C$) to the corner square.
\end{lemma}

\proof. Assume that there exists a valid $n\times n$ ($n > 3$)
board configuration with a $3\times 3$ corner block $C$ that
contains only two guards.  Without loss of generality, suppose that
$C$ is situated in the upper left corner of the board.  Since four
guards are required to protect a prisoner residing on an interior
square, and we have two guards, a guard must occupy $c_{22}$. Since,
by assumption, there remains only one more guard, there must lie a
prisoner in   $c_{12}$ or $c_{21}$.  However, three guards are
required to cover a prisoner that is situated in a non-corner edge
square. Hence, the board cannot be valid and we have a
contradiction.

\begin{figure}[ht]
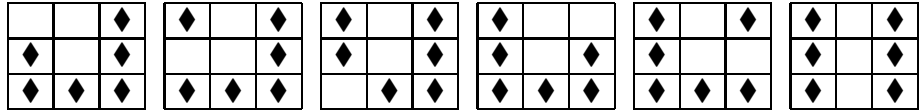

\begin{tabular}{|r|r|r|}
\hline
\ & \ \  & \b \\
\hline
\b & \ \  & \b \\
\hline
\b &\b \  & \b \\
\hline
\end{tabular}\ \
\begin{tabular}{|r|r|r|}
\hline
\b & \ \  & \b \\
\hline
\ & \ \  & \b \\
\hline
\b &\b \  & \b \\
\hline
\end{tabular}\ \
\begin{tabular}{|r|r|r|}
\hline
\b & \ \  & \b \\
\hline
\b & \ \  & \b \\
\hline
\ &\b \  & \b \\
\hline
\end{tabular}\ \
\begin{tabular}{|c|c|c|}
\hline
\b  & \ \  & \ \\
\hline
\b & \  & \b \\
\hline
\b &\b & \b \\
\hline
\end{tabular}\ \
\begin{tabular}{|c|c|c|}
\hline
\b  & \ & \b  \\
\hline
\b & \   & \  \\
\hline
\b &\b & \b \\
\hline
\end{tabular}\ \
\begin{tabular}{|r|r|r|}
\hline
\b  & \ \  & \b \\
\hline
\b & \ \  & \b \\
\hline
\b &\ \  & \b \\
\hline
\end{tabular}
\caption{Possible $3\times 3$ UL Corner Blocks With 3 Guards}
\label{possiblecornerblocks}
\end{figure}

To establish the last part of the claim, observe that besides the
guard at $c_{22}$ the other two are either at $c_{12}$ or $c_{21}$,
or they must lie adjacent to the prisoner located at $c_{12}$ or
$c_{21}$.  (If one of the guards were located in $c_{33}$ then we
would have a prisoner in either $c_{12}$ or $c_{21}$ without a
sufficient number of adjacent guards).  Possible arrangements, up to
a reflection about the main diagonal, appear in
Figure~\ref{possiblecornerblocks}.  In these cases, we are left with
a prisoner in $c_{33}$.  \eproof

\begin{figure}[ht]
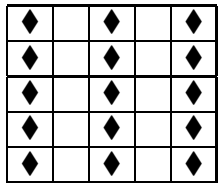

\begin{tabular}{|r|r|r|r|r|}
\hline
\b  \       & \  \  & \b  & \ \ & \b  \\
\hline
\b  \       & \  \  & \b  & \ \ & \b  \\
\hline
\b  \       & \  \  & \b  & \ \ & \b  \\
\hline
\b  \       & \  \  & \b  & \ \ & \b  \\
\hline
\b  \       & \ \   & \b  & \ \ & \b  \\
\hline
\end{tabular}
\caption{The Maximum $5\times 5$ Board Configuration}
\label{5x5maxboard}
\end{figure}

Now, we are ready to consider the $5\times 5$ case.  The only
maximum configuration (up to rotations) is illustrated in
Figure~\ref{5x5maxboard}.  Proposition~\ref{p5} establishes that
this is a maximum $5\times 5$ board configuration.

\begin{proposition}\label{p5}
A maximum $5\times 5$ board configuration contains fifteen
prisoners; that is, $P(5)=15$.
\end{proposition}

\n \proof.  Assume that there exists a valid $5\times 5$ board
configuration with 16 or more prisoners.  We will see that this
leads to a contradiction.

\begin{figure}[h]
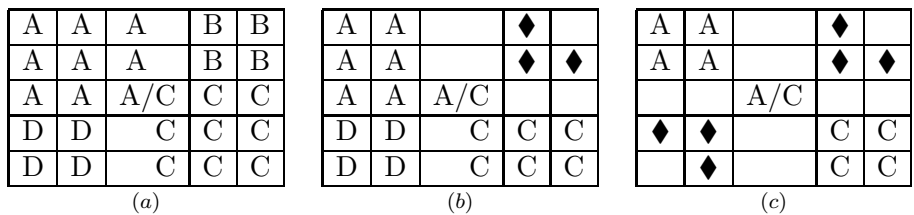

\[\underset{(a)}{ \begin{tabular}{|c|c|r|c|c|} \hline
A & A & A \ \ \   & B & B \\
\hline
A & A & A \ \ \   & B & B   \\
\hline
A & A & A/C  & C & C  \\
\hline
D & D & \ C   & C & C   \\
\hline
D & D & \ C & C & C  \\
\hline
\end{tabular}}\ \ \ \
\underset{(b)}{
\begin{tabular}{|c|c|r|c|c|} \hline
A & A & \ & \b & \ \\
\hline
A & A & \ & \b & \b \\
\hline
A & A & A/C  & \ & \  \\
\hline
D & D & \ C   & C & C   \\
\hline
D & D & \ C & C & C  \\
\hline
\end{tabular}}\ \ \ \
\underset{(c)}{
\begin{tabular}{|c|c|r|c|c|} \hline
A & A & \ & \b & \ \\
\hline
A & A & \ & \b & \b \\
\hline
\ & \ & A/C  & \ & \  \\
\hline
\b & \b & \ & C & C   \\
\hline
\ & \b & \ & C & C  \\
\hline
\end{tabular}}
\]
\caption{Partitions of the $5\times 5$ Board}
\label{5x5boards}
\end{figure}

Divide the $5\times 5$ board into two opposite (overlapping) corner
$3\times 3$ blocks, $A$ and $C$, that have a square in common and two $2\times 2$
opposite corner blocks, $B$ and $D$, as illustrated in Figure~\ref{5x5boards}(a).
\medskip

According to Lemma~\ref{threebythree}, the two $3\times 3$ blocks
$A$ and $C$ collectively contain at most $2(6)-1=11$ prisoners since
the shared square (common to blocks $A$ and $C$) must contain a
prisoner.  Recall that Lemma~\ref{twobytwo} establishes that each of
blocks $B$ and $D$ contains at most three prisoners.  Thus, for the
board to contain a total of sixteen prisoners we must find either
ten or eleven prisoners shared in blocks $A$ and $C$.  (They cannot
share just nine prisoners since that would force a $2\times 2$ block
to hold four prisoners).

{\bf Case 1.  Blocks $A$ and $C$ share 11 prisoners.} In this case,
one of the $2\times 2$ blocks holds three prisoners and the other
holds two prisoners.  Without loss of generality, let us suppose
that the $B$ block has three prisoners; then the one guard must lie
in the (1,5) position. To cover the three prisoners in block $B$,
guards in the (1,3), (2,3), (3,4), and (3,5) positions (refer to
Figure~\ref{5x5boards}(b)).  For the board to be valid, block $A$
must match one of the corner blocks depicted in
Figure~\ref{possiblecornerblocks}, none of which allows guards in
both the (1,3) and the (2,3) positions.  Therefore, the board is not
valid and we have a contradiction in this case.

{\bf Case 2.  Blocks $A$ and $C$ share 10 prisoners.} In this case,
each of the $2\times 2$ blocks $B$ and $D$ holds three prisoners. In
order to maintain a valid board configuration, we are then forced to
place guards in the (1,3), (1,5), (2,3), (3,1), (3,2), (3,4), (3,5),
(4,3), and (5,3) positions as indicated in
Figure~\ref{5x5boards}(c).  But then there remain only nine
uncommitted squares in which to place the ten prisoners that blocks
$A$ and $C$ are supposed to share.  Thus, we also find a
contradiction in this case.  \eproof
\bigskip

\begin{figure}[ht]
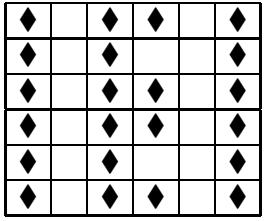

\begin{tabular}{|r|r|r|r|r|r|}
\hline
\b & \ \  & \b  & \b & \ \ & \b \\
\hline
\b & \ \   & \b  & \ \ &\ \  &\b      \\
\hline
\b & \ \   & \b   & \b & \ \ &\b \\
\hline
\b & \ \   & \b  & \b &\ \ &\b  \\
\hline
\b & \ \   & \b  & \ \ &\ \ &\b \\
\hline
\b  & \ \  & \b  & \b &\ \ &\b  \\
\hline
\end{tabular}
\caption{A Maximum $6\times 6$ Board Configuration}
\label{6x6maximumboard}
\end{figure}

For $n=6$ all maximum boards amount to rotations or small
perturbations of the arrangement in Figure~\ref{6x6maximumboard},
the validity of which yields the lower bound $P(6)\ge 22$. We will
show that, in fact, $P(6)=22$ by an analysis of manageable size.
Dunbar, Hoffman, Laskar, and Markus assert (without proof) a fact
about 1/2-domination in the king's graph dimension $6\time 6$ which,
if true, implies that $P(6)=22$ (see \cite{dhlm}). This is indeed
the case, as we shall establish next. We use a more specific version
of Lemma~\ref{threebythree} in order to obtain this fact.

\begin{lemma}\label{threebythreees}
If $C$ is a $3\times 3$ corner block holding six prisoners within a
valid board ($n > 3$), where $c_{11}$ is the corner square, then
up to a diagonal symmetry the block has one of the six arrangements
in Figure~\ref{possiblecornerblocks}.
\end{lemma}

\n \proof.  The position $c_{22}$ must have a guard as we have seen
and also one of the positions $c_{12}$ or $c_{21}$ must have a
guard. By symmetry we can assume we have a guard at $c_{12}$.  Then
there are six possible spots for the third guard. This gives exactly
the six arrangements in Figure~\ref{possiblecornerblocks}. \eproof

\begin{proposition}\label{p6}
A maximum $6\times 6$ board contains twenty two prisoners.  That is, $P(6)=22$.
\end{proposition}

\begin{figure}[ht]
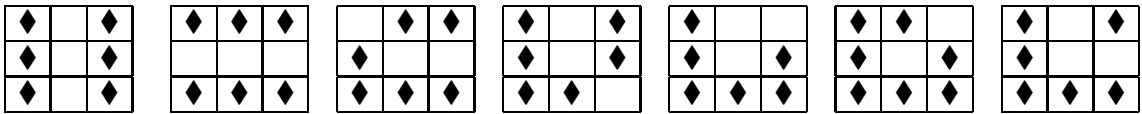

\begin{tabular}{|r|r|r|}
\hline
\b  & \ \  & \b \\
\hline
\b & \ \  & \b \\
\hline
\b &\ \  & \b \\
\hline
\end{tabular} \ \ \
\begin{tabular}{|r|r|r|}
\hline
\b & \b \  & \b \\
\hline
\ & \ \  & \ \\
\hline
\b &\b \  & \b \\
\hline
\end{tabular}\ \ \
\begin{tabular}{|r|r|r|}
\hline
\ & \b \  & \b \\
\hline
\b & \ \  & \ \\
\hline
\b &\b \  & \b \\
\hline
\end{tabular}\ \ \
\begin{tabular}{|r|r|r|}
\hline
\b & \ \  & \b \\
\hline
\b & \ \  & \b \\
\hline
\b &\b \  & \ \\
\hline
\end{tabular}\ \ \
\begin{tabular}{|r|r|r|}
\hline
\b & \ \  & \ \\
\hline
\b & \ \  & \b \\
\hline
\b &\b \  & \b \\
\hline
\end{tabular}\ \ \
\begin{tabular}{|r|r|r|}
\hline
\b & \b \  & \ \\
\hline
\b & \ \  & \b \\
\hline
\b &\b \  & \b \\
\hline
\end{tabular}\ \ \
\begin{tabular}{|r|r|r|}
\hline
\b & \ \  & \b \\
\hline
\b & \ \  & \ \\
\hline
\b &\b \  & \b \\
\hline
\end{tabular}
\caption{Possible Upper Right $3\times 3$ Corner Blocks}
\label{corneroptions}
\end{figure}

\n \proof. We have observed that $P(6)\ge 22$. To verify that
$P(6)\le 22$ let us assume the existence of a valid arrangement $C$
with 23 prisoners; we shall see that this leads to a contradiction.
By Lemma~\ref{threebythree}, three of the four $3\times 3$ corner
blocks have six prisoners and one has five prisoners. Without loss
of generality, we may assume that the block with five prisoners is
the lower right one. By Lemma~\ref{threebythreees} and by symmetry,
we can assume that the block in the upper left corner is one of
those in Figure~\ref{possiblecornerblocks} and the upper right
$3\times 3$ corner block is an arrangement found in
Figure~\ref{corneroptions}. Since one of these arrangements  has
diagonal symmetry we have only seven possible situations as listed
Figure~\ref{corneroptions}.

In each case one can find that the prisoner at $c_{13}$ or $c_{24}$
does not have enough guards around it. This contradicts the
existence of a configuration with 23 prisoners. \eproof

We suspect that this block partition approach can be adapted to
compute or bound $P(n)$ for larger sizes of $n$, although this
approach could turn out to be quite lengthy. These proofs may very
well be pursued as undergraduate research projects.

\section{\label{deficiency}Upper bound for $P(n)$ and the deficiency function}

As the board size grows larger, establishing the exact number of
prisoners on a maximum board becomes increasingly difficult.  The
proof of Proposition~\ref{maximum4x4boards} foreshadows the
importance of finding useful upper bounds on $P(n)$. In this
section, we construct a tool that will help in establishing these
bounds -- the deficiency matrix.  We then use the deficiency matrix
to determine a general upper bound for $P(n)$.

Suppose that we have fixed the board size at $n\times n$ ($n \ge
3$). With each configuration, we associate a binary matrix
$X=\left(x_{ij}\right)$ defined by

$$x_{ij}=\left\{
\begin{array}{cl}
1 & \mbox{if a prisoner lies in the $(i,j)$ position}\\
0 & \mbox{if a guard lies in the $(i,j)$ position.}
\end{array}
\right.$$

\noindent  Many who work in combinatorics and graph theory, such as
Hedetniemi, Hedetniemi, and Reynolds (see \cite{hedetniemi}) have
employed this idea. In any local measure of optimality we must be
attentive to the number of prisoners lying in square adjacent to a
particular square; we let $x_{ij}^*$ denote the number of prisoners
lying in squares adjacent to the $(i,j)$ square.

The deficiency matrix serves as an ad-hoc, local measure of the
optimality of a given board configuration.  Its construction arises
from our observations and conjectures of maximum board
configurations.  We define the deficiency matrix
$\delta=\left(\delta_{ij}\right)$ by
$$\delta_{ij}=\mbox{expectation}-x_{ij}^*, \mbox{where}$$
$$\mbox{expectation}=\left\{
\begin{array}{cll}
1 & \mbox{if $(i,j)$ is a corner square with} & x_{ij}=1\\
2 & \mbox{if $(i,j)$ is a corner square with} & x_{ij}=0\\
2 & \mbox{if $(i,j)$ is an edge square with} & x_{ij}=1\\
4 & \mbox{if $(i,j)$ is an edge square with} & x_{ij}=0\\
4 & \mbox{if $(i,j)$ is an interior square with} & x_{ij}=1\\
6 & \mbox{if $(i,j)$ is an interior square with} & x_{ij}=0.
\end{array}\right.$$

Figure~\ref{4x4deficiency_example} depicts a $4\times 4$ non-maximal
board configuration and its corresponding deficiency matrix.  The
positive entries in the deficiency matrix indicate areas of the
board that are thought to be less than optimal; the 2's indicate
that the ``worst'' deficiencies occur on the corresponding interior
squares.
\bigskip

\begin{figure}[h]
\begin{tabular}{|c|c|c|c|}
\hline
\b & \   & \b & \ \\
\hline
\ & \b & \ & \b \\
\hline
\b & \   & \b & \ \\
\hline
\ & \b & \ & \b \\
\hline
\end{tabular} \  \  \
\begin{tabular}{|c|c|c|c|}
\hline
$0$ & $1$ & $0$ & $0$ \\
\hline
$1$ & $0$ & $2$ & $0$ \\
\hline
$0$ & $2$ & $0$ & $1$ \\
\hline
$0$ & $0$ & $1$ & $0$ \\
\hline
\end{tabular}
\caption{Non-maximal $4\times 4$ Board and Its Deficiency Matrix}
\label{4x4deficiency_example}
\end{figure}
\bigskip

Since we use these values to obtain an upper bound on $P(n)$, it
helps to first consider bounds on $\delta_{ij}$ for $1\le i,j\le n$.
Suppose that the $(i,j)$ square contains a prisoner.  If $(i,j)$ is
a corner square, then at least two of the three adjacent squares
must contain guards; therefore $x_{ij}^*\le 1$ and in checking our
expectation value above we see that $\delta_{ij}\ge 0$.  Likewise,
if $(i,j)$ is an edge square containing a prisoner or an interior
square with a prisoner, we find that $\delta_{ij}\ge 0$.

Suppose that we find a guard in a corner square $(i,j)$.  Then three
squares lie adjacent to this square, so we find at most three
prisoners in the neighboring squares.  Therefore, $x_{ij}^*\le 3$
and so $\delta_{ij}\ge -1$.  Via similar considerations, we find
that if a guard occupies an edge square then $\delta_{ij}\ge -1$ and
for an interior square we get $\delta_{ij}\ge -2$.

We define the {\it net deficiency of a board configuration} as the
sum of all entries in the deficiency matrix,
$$\Delta=\sum_{i,j=1}^n\delta_{ij}.$$ \noindent For instance, the
board configuration in Figure~\ref{4x4deficiency_example} has a net
deficiency of 8.  We expect maximum board configurations to
correspond to minimum net deficiency values.  We will relate
$\Delta$ to the overall number of guards in a given board
configuration.  Let $P_C$ and $G_C$ denote the total number of
prisoners and guards, respectively, found in the corner squares.
Similarly, $P_E$ and $G_E$ refer to the prisoners and guards in edge
squares, and $P_I$ and $G_I$ refer to prisoners and guards in
interior squares.  With this notation we have
$$\begin{array}{rll}
\Delta & = & \sum_{\mbox{corners}}\delta_{ij}+
\sum_{\mbox{edges}}\delta_{ij}+
\sum_{\mbox{interior}}\delta_{ij} \\
\ & \ge & \left(0\cdot P_C-1\cdot G_C\right) +\left(0\cdot P_E-1\cdot G_E\right)+\left(0\cdot P_I - 2\cdot G_I\right)\\
\  & = & -G_C-G_E-2\cdot G_I.
  \end{array}
$$

\noindent This establishes the next lemma.

\begin{lemma}  The net deficiency of a given configuration satisfies the inequality $\Delta \ge -G_C-G_E-2\cdot G_I$.\label{deficiencybound}\end{lemma}

Now we are ready to think about bounding the size of $P(n)$.  Observe that
\begin{equation}\label{ineeq}
\begin{array}{l} 4x_{ij}+x_{ij}^{*}=
\begin{cases} 8-\delta_{ij}, \mbox{if $x_{ij}=1$  and $(i,j)$ is an interior square } \\
6-\delta_{ij}, \mbox{if $x_{ij}=0$ and $(i,j)$ is an interior square}
\end{cases}\\  \\ 3x_{ij}+x_{ij}^{*}=
\begin{cases} 5-\delta_{ij}, \mbox{ if $x_{ij}=1$ and $(i,j)$ is an edge square} \\
 4-\delta_{ij}, \mbox{if $x_{ij}=0$ and $(i,j)$ is an edge square}
\end{cases}
\\  \\ 2x_{ij}+x_{ij}^{*}=
\begin{cases} 3-\delta_{ij}, \mbox{if $x_{ij}=1$ and $(i,j)$ is a corner square} \\
 2-\delta_{ij}, \mbox{if $x_{ij}=0$ and $(i,j)$ is a corner square.}
\end{cases}
\end{array}
\end{equation}

\begin{theorem}\label{main}  The number of prisoners in a valid configuration is given by
\begin{equation}\label{exacteqnk}
P=\frac{3n^2}{5}-\frac{4n}{5}+\frac{1}{10}(3P_E+6P_C-\Delta ).
\end{equation}
\end{theorem}

\proof.  \ \ Summing the left hand sides of the equations in
(\ref{ineeq}) over all squares of the board, we obtain $$4\cdot P_I
+3\cdot P_E+2\cdot P_C+\sum_{1\le i,j\le n}x_{ij}^*=4\cdot P_I
+3\cdot P_E+2\cdot P_C+8\cdot P_I +5\cdot P_E+3\cdot P_C=12\cdot P_I
+8\cdot P_E+5\cdot P_C.$$

We will equate this result with the sum of the right hand sides. In
summing over the interior squares that contain prisoners, this
contributes $8-\delta_{ij}=6+2-\delta_{ij}$ for each such square ,
whereas the interior square that contain guards contribute only
$6-\delta_{ij}$ per square.  There are $(n-2)^2$ interior squares,
so altogether these sum to $6(n-2)^2+2\cdot P_I-\sum_{\text{int.
sqrs.}}\delta_{ij}$.  Similarly summing the right hand sides over
all edge squares we get $4\cdot 4\left[4(n-2)\right]+1\cdot P_E-
\sum_{\text{edge sqrs.}}\delta_{ij}$.  Summing over the corners
yields $8+1\cdot P_C-\sum_{\text{corner sqrs.}}\delta_{ij}$.
Combining these right-hand sums and equating with the left-hand sum,
we obtain the equation

$$12P_I+8P_E+5P_C=6(n-2)^2+2P_I+16(n-2)+P_E+8+P_C-\Delta $$
or
$$10P_I+7P_E+4P_C=6n^2-8n-\Delta.$$
Then since $P=P_I+P_E+P_C$ we then obtain

$$10P=6n^2-8n+3P_E+6P_C-\Delta,$$
which leads to (\ref{exacteqnk}).\eproof

By combining the inequality in Lemma~\ref{deficiencybound} with this
theorem, we obtain a crude upper bound on $P(n)$.

\begin{corollary}\label{firstupb} In a maximum configuration of prisoners and guards on a $n\times n$ board the number of prisoners obeys the inequality
\begin{equation}\label{upperbfkg}
P(n) \le \frac{2n^2+n}{3}.
\end{equation}
\end{corollary}

\proof. \  Using Lemma~\ref{deficiencybound} and (\ref{exacteqnk}) we obtain
\begin{eqnarray*}
\ds P(n) & \le & \frac{3n^2}{5}-\frac{4n}{5}+\frac{1}{10}(3P_E+6P_C+2G_I+G_E+G_C)\\
& = & \ds \frac{3n^2}{5} - \frac{4n}{5} + \frac{1}{10}\Big[2P_E+5P_C+2(n-2)^2-2P_I+4(n-2)+4\Big]\\
& = & \ds \frac{3n^2}{5} - \frac{4n}{5} + \frac{1}{10}\Big[2P_E+5P_C+2(n-2)^2-2P+2P_E+2P_C+4n-4\Big].
\end{eqnarray*}

\n This implies
\begin{eqnarray*}
\ds \left(1+\frac{1}{5}\right )P(n) & \le & \ds \frac{3n^2}{5} - \frac{4n}{5} + \frac{1}{10}(4P_E+7P_C+2n^2-4n+4)\\
& \le & \ds
\frac{3n^2}{5}-\frac{4n}{5}+\frac{1}{10}\Big[(4)(4)(n-2)+(7)(4)+2n^2-4n+4\Big]\\
& = & \frac{4n^2+2n}{5}
\end{eqnarray*}

\n Therefore $\ds P(n)\le \left(\frac{5}{6}\right)\left(\frac{4n^2+2n}{5}\right)=\frac{2n^2+n}{3}$.\eproof

We believe that $\Delta \le O(n)$ in general.  This fact is
equivalent to $P(n) \le 3n^2/5 +O(n)$. However we can tighten the
upper bound (\ref{upperbfkg}) by getting a better estimate for
$\Delta$.

\begin{lemma}\label{secondub} In a valid configuration the net deficiency satisfies
$$\Delta\ge -1\cdot G = -1(G_I+G_E+G_C).$$
\end{lemma}

\n \proof.  \ Recall our previous observations about the possible
range of values for $\delta_{ij}$.  If the $(i,j)$ board position
contains a prisoner then from the definition it follows that
$\delta_{ij}\ge 0$.  If the square is a corner or edge square
containing a guard then $\delta_{ij}\ge -1$.  For an interior square
containing a guard, we have noted that $\delta_{ij}\ge -2$.  Let us
focus on this last case.

Suppose that a guard occupies the $(i,j)$ interior position in a
valid board configuration and that $\delta_{ij}=-2$.  Then it must
be the case that all adjacent squares contain prisoners, as depicted
in Figure~\ref{configurationnearguard}(a).  The g's denote guards
that are then forced into the arrangement in order for the
configuration to be valid.  We see that each of the prisoners in the
squares diagonally adjacent to this position lies adjacent to five
or six guards (depending on the occupants in the squares marked with
asterisks).  The possible deficiency values for neighboring squares
appear in Figure~\ref{configurationnearguard}(b).  Summing these
deficiency values, we find that the net contribution of the $3\times
3$ block satisfies $2\le \Delta_{local} \le 6$.  Notice that, as the
g's in Figure~\ref{configurationnearguard}(a) suggest, it is not
possible for two such $3\times 3$ blocks around guards with
deficiency $-2$ to overlap. Thus, we see that each guard on the
board contributes a net deficiency value not less than $-1$.

\begin{figure}[ht]
\[\underset{(a)}{
\begin{tabular}{|c|c|c|c|c|}
\hline
$\ast$ & g & g & g & $\ast$  \\
\hline
g & \b  & \b \ & \b &g \\
\hline
g & \b  & \ \ & \b &g \\
\hline
g & \b  & \b \ & \b &g \\
\hline
$\ast$ & g & g & g & $\ast$  \\
\hline
\end{tabular}}
\ \ \ \
\underset{(b)}{
\begin{tabular}{|c|c|c|c|c|}
\hline
\ & \ & \ & \ & \  \\
\hline
\ \ \ & {\Small 1, 2} & 0 & {\Small 1, 2} & \ \ \ \\
\hline
\ & 0 & -2 & 0 & \ \\
\hline
\ & {\Small 1, 2} & 0 & {\Small 1, 2} & \ \\
\hline
\ & \ & \ & \ & \  \\
\hline
\end{tabular}}
\]
\caption{Local Configuration Near a Guard with Deficiency -2}
\label{configurationnearguard}
\end{figure}

Summing the $\delta_{ij}$'s over all board positions, we have
$\Delta = \sum_{prisoners}\delta_{ij} + \sum_{guards}\delta_{ij}
\ge -1\cdot G$.  \eproof

Using this bound on $\Delta$ in Theorem~\ref{main}, we obtain a
better upper bound for $P(n)$.  The calculations parallel those used
in the proof of Corollary~\ref{firstupb}.

\begin{theorem}\label{secondubforp} For an $n\times n$ maximum arrangement of prisoners
and guards, the number of prisoners, $P(n)$, satisfies the inequality
\begin{equation}\label{mainineq}
P(n)\le \frac{7n^2+4n}{11}.
\end{equation}
\end{theorem}

\n \proof.  By Lemma~\ref{secondub}, $-\Delta\le G$.  Applying this
upper bound in Theorem~\ref{main}, we get
\begin{eqnarray*}
P & \le & \frac{3}{5}n^2 - \frac{4}{5}n + \frac{1}{10}\Big[3P_E+6P_C+G \Big] = \frac{3}{5}n^2 - \frac{4}{5}n + \frac{1}{10}\Big[3P + 3P_C-3P_I+G \Big]\\
& = & \frac{3}{5}n^2 - \frac{4}{5}n + \frac{1}{10}\Big[2P+3\big(P_C-P_I\big)+n^2 \Big].
\end{eqnarray*}
\n Subtracting $\frac{2}{10}P$ from both sides and combining the $n^2$ terms, we see that this implies
\begin{eqnarray*}
\frac{8}{10}P &\le & \frac{7}{10}n^2 - \frac{4}{5}n + \frac{3}{10}\big(P_C-P_I\big)\\
&\le &  \frac{7}{10}n^2 - \frac{4}{5}n + \frac{3}{10}\big(4-P+P_E+P_C\big)\\
&\le & \frac{7}{10}n^2 - \frac{4}{5}n + \frac{3}{10}\big(4-P+4n-4\big).
\end{eqnarray*}
\n The claim now follows after a bit of arithmetic.\eproof

It seems that this method of finding an upper bound can be further
sharpened, but we have not found a complete analysis to show that
$\Delta \ge -O(n)$. This would be the sharpest result for this type
of domination problem.

\section{\label{whatelse}Other Results, Conjectures and Open Questions}

Ionascu, Pritikin, and Wright have established values of $P(n)$ for
$n\in \{7, 8, 9, 10, 11\}$ (see \cite{ipw}).  Most of their
arrangements were obtained using the LPSolve IDE program in the
Lesser GNU public domain for solving integer linear programming
problems with Branch-and-Bound and Simplex Methods. The second
author used CPLEX while visiting at the Georgia Institute of
Technology in the Faculty Development Program in 2005-2006; with the
help of Professor William Cook he analyzed the case $n=11$.

Many interesting questions remain to be answered.  What are the
values of $P(n)$ for integers $n$ larger than 11?  With error-free
play, does one particular player enjoy an advantage?  Perhaps the
advantage varies with the board size.

If $P(n)$ is odd, we conjecture that the game favors the red player,
but it is not clear that a winning strategy exists. When $P(n)$ is
even, we suspect that error-free play by both players will lead to a
tie.

Given that we find several maximal $4\times 4$ board configurations
with eight prisoners (an even number), it seems that the second
player (blue) will find opportunities to win unless s/he is forced
to use Rule I. The question is: can the red player always achieve a
win or a tie? We believe there is a strategy for the red player to
win despite all of these chances for the blue player.  In general,
it is apparent that that the final maximal configuration is an
important factor in the game, since the number of prisoners in it
determines the fate of the game. So it is in the red player's
interest to end in a maximal arrangement with an odd number of
prisoners on the board. Similarly it is part of blue strategy to
divert the end configuration to a maximal one that has an even
number of prisoners. Each player can change the configuration at
only one place at which the opponent has already directed the game
toward his final configuration and leave one place as it is. As a
result, almost half the prisoners on the final board configuration
are where each player wanted them to be. So from this perspective
the end game is dictated by the parity and the number of maximal
configurations with $P(n)-1$, $P(n)-2$, ... prisoners.

In this paper, we have shown that $P(n)$ is bounded above by
$\ds\frac{7n^2+4n}{11}$, but we conjecture that $P(n)=3n^2/5+O(n)$.
It would be great to see someone, especially a student, improve our
upper bound.

{\bf Acknowledgments. } Special thanks go to Annalisa Crannell, of
Franklin and Marshall College, for her many helpful comments and
suggestions for this paper. We also thank  Mike McCoy, a Columbus
State Univeresity student who spent hours designing a program that
has a very good run-time, even for large size boards, and is
producing arrangements that we think are ``very close" to  maximum
arrangements.

\end{document}